\documentclass[11pt]{amsart}
\usepackage{amssymb}

\usepackage{graphicx}
\usepackage{amscd}

\newtheorem{theorem}{Theorem}
\theoremstyle{plain}

\newtheorem{definition}{Definition}

\newtheorem{proposition}{Proposition}
\newtheorem{remark}{Remark}

\numberwithin{equation}{section}

\begin{document}
\title[Generalized functions and differential operators]{Colombeau generalized functions and solvability of differential operators}
\author{Khaled Benmeriem}
\address{Centre Universitaire de Mascara. Mascara, Algeria.}
\email{benmeriemkhaled@yahoo.fr}
\author{Chikh Bouzar}
\address{Department of Mathematics. University of Oran Essenia. Algeria.}
\email{bouzar@wissal.dz}
\maketitle

\begin{abstract}
The aim of this paper is to prove that the well known non solvable Mizohata
type partial differential equations have Colombeau generalized solutions
which are distributions if and only if they are solvable in the space of
Schwartz distributions. Therefore the Colombeau generalized solvability
includes both a new solution concept and new mathematical objects as
solutions.
\end{abstract}

\section{Introduction}

Colombeau generalized functions were introduced , see \cite{col}, in
connection with the so-called problem of multiplication of Schwartz
distributions \cite{Schw}. They were developed and applied in important
nonlinear problems, see \cite{Biag}, \cite{col2} and \cite{o}. General
methods of construction of such generalized functions were given in \cite
{A-R} \ and \cite{M}. The authors of \cite{NPS} have tackled the linear
counterpart of this theory.

The theory of Colombeau generalized functions provides new solutions of
partial differential equations; these new solutions can be divided into two
categories:

1) there are classical functions or distributions which are solutions (in
one of the new senses provided by this theory) of partial differential
equations without solution in the sense of distributions, e. g. see \cite
{Biag}, \cite{col2}, \cite{Hu}, \cite{col-ober} and \cite{col-he}.

2) there are also new objects (such as the square of the Dirac delta
distribution,...) which can be solutions of equations.

In \cite{col-he} the fundamental concept of regularized derivatives was
studied and results on global solvability, in the framework of this theory,
of the Cauchy problem for large classes of regularized partial differential
equations have been given. In particular, the well-known non solvable
Mizohata differential equations with regularized derivatives become solvable
in the Colombeau algebra. It is then interesting to show the relation
between Colombeau generalized solutions and distributional solutions if they
exist.

The paper deals, in the framework of the simplified Colombeau algebra, with
a class of differential operators non solvable in distributions theory. We
show that their Colombeau generalized solutions as regularized differential
equations are in relations with distributional solutions if and only if they
are solvable in the space of Schwartz distributions. Therefore in the
general case in which there are no distributional solution, the new
solutions from \cite{col-he} are not associated with classical objects, even
if they are solutions in a new sense: an enlargement of the reservoir of
mathematical objects that could be solutions is really needed.

\textbf{Acknowledgements:} The authors thank Professor J.-F.
Colombeau for the suggested amelioration of the paper and
Professor M. Oberguggenberger for the useful comments on Theorem
4.

\section{Simplified algebra of Colombeau \label{par1}}

In this section we recall the simplified Colombeau algebra of generalized
functions and some needed notions of this theory, for a deep study see \cite
{col}, \cite{col2} and \cite{o}. Let $\Omega $ be a non void open subset of $%
\mathbb{R}^{d}$ and $I=$ $\left] 0,1\right[ ,$ define $\chi _{M}\left(
\Omega \right) $ as the space of elements $\left( u_{\varepsilon }\right)
_{\varepsilon }$ of $\chi \left( \Omega \right) =\left( C^{\infty }\left(
\Omega \right) \right) ^{I}$ such that, for every compact set $K\subset
\Omega $, $\forall \alpha \in \mathbb{Z}_{+}^{d}$, $\exists m>0,$
\begin{equation}
\,\underset{x\in K}{\sup }\left| \partial ^{\alpha }u_{\varepsilon }\right|
\leq O\left( \varepsilon ^{-m}\right) ,\text{ as }\varepsilon \rightarrow 0.
\label{2.1}
\end{equation}
By $\mathcal{N}\left( \Omega \right) $ we denote the elements $\left(
u_{\varepsilon }\right) _{\varepsilon }\in \chi _{M}\left( \Omega \right) $
satisfying $\forall K\subset \Omega $, $\forall \alpha \in \mathbb{Z}%
_{+}^{d} $, $\forall q>0,$
\begin{equation}
\,\underset{x\in K}{\sup }\left| \partial ^{\alpha }u_{\varepsilon }\right|
\leq O\left( \varepsilon ^{q}\right) ,\text{ as }\varepsilon \rightarrow 0.%
\text{ }  \label{2.2}
\end{equation}
An element of $\chi _{M}\left( \Omega \right) $ is called moderate and an
element of $\mathcal{N}\left( \Omega \right) $ is called null. It is easy to
prove that $\chi _{M}\left( \Omega \right) $ is an algebra and $\mathcal{N}%
\left( \Omega \right) $ is an ideal of $\chi _{M}\left( \Omega \right) .$

\begin{definition}
The simplified algebra of Colombeau defined on $\Omega $, denoted $\mathcal{G%
}_{s}\left( \Omega \right) ,$ is the quotient algebra
\begin{equation*}
\mathcal{G}_{s}\left( \Omega \right) =\frac{\chi _{M}\left( \Omega \right) }{%
\mathcal{N}\left( \Omega \right) }.
\end{equation*}
\end{definition}

The algebra of Colombeau $\mathcal{G}_{s}\left( \Omega \right) $ is a
commutative and associative differential algebra containing $D^{\prime
}\left( \Omega \right) $ as a subspace and $C^{\infty }\left( \Omega \right)
$ as subalgebra, see for details \cite{col}, \cite{col2} and \cite{o}, where
others important properties of this algebra are studied.

Recall the notion of association relation in the Colombeau algebra $\mathcal{%
G}_{s}\left( \Omega \right) $, a generalized function $u\in \mathcal{G}%
_{s}\left( \Omega \right) $ and a distribution $T$ $\in D^{\prime }\left(
\Omega \right) $ are called associated, denoted $u\approx T$, if there
exists $\left( u_{\varepsilon }\right) _{\varepsilon }$ a representative of $%
u$ such that, $\forall \phi \in C_{0}^{\infty }\left( \Omega \right) ,$
\begin{equation*}
\underset{\varepsilon \rightarrow 0}{\lim }\int u_{\varepsilon }\left(
x\right) \phi \left( x\right) =\left\langle T,\phi \right\rangle .
\end{equation*}

We introduce, for our need, an association relation less stronger than the
classical association.

\begin{definition}
A generalized function $u\in \mathcal{G}_{s}\left( \Omega \right) $ and a
distribution $T$ $\in D^{\prime }\left( \Omega \right) $ are called locally
associated at $x_{0}\in \Omega $, denoted $u\approx _{x_{0}}T$, if there
exists $\left( u_{\varepsilon }\right) _{\varepsilon }$ a representative of $%
u $ and $\omega \subset \Omega $ an open neighborhood of $x_{0},$ such that,
$\forall \phi \in C_{0}^{\infty }\left( \omega \right) ,$
\begin{equation*}
\underset{\varepsilon \rightarrow 0}{\lim }\int u_{\varepsilon }\left(
x\right) \phi \left( x\right) =\left\langle T,\phi \right\rangle .
\end{equation*}
\end{definition}

The proof of the following result is easy.

\begin{proposition}
Let $u\in \mathcal{G}_{s}\left( {\Omega }\right) $ and $T\in D^{\prime
}\left( {\Omega }\right) ,$ then $u\approx _{x_{0}}T,\forall x_{0}\in {%
\Omega },$ if and only if $u\approx T.$
\end{proposition}

\section{Regularized partial differential equations}

For the concept of regularized derivatives of Colombeau generalized
functions and its application to general Cauchy problems see \cite{col-he}.
Denote by $\mathcal{H}$ the set of non-decreasing functions $h:I\rightarrow
I $, such that $\lim\limits_{\varepsilon \rightarrow 0}h\left( \varepsilon
\right) =0.$ Let $\rho \in C_{0}^{\infty }\left( \mathbb{R}^{d}\right) $ and
$\int \rho \left( x\right) dx=1,$ we define the sequence $\left( \rho
_{\varepsilon }\right) _{\varepsilon }$ by $\rho _{\varepsilon }\left(
x\right) =\dfrac{1}{\varepsilon ^{d}}\rho \left( \dfrac{x}{\varepsilon }%
\right) ,\varepsilon \in I.$

\begin{definition}
Let $u\in \mathcal{G}_{s}\left( \mathbb{R}^{d}\right) $ and $h\in \mathcal{H}
$, the partial regularized derivative of $u$ with respect to $x_{j}$,
denoted $\left( \widetilde{\partial }_{x_{j}}\right) _{h}u,$ is defined by
\begin{equation*}
\left( \widetilde{\partial }_{x_{j}}\right) _{h}u=cl\left( \partial
_{x_{j}}u_{\varepsilon }\ast \rho _{h\left( \varepsilon \right) }\right)
_{\varepsilon \in I}\text{ \ ,}
\end{equation*}
where $\left( u_{\varepsilon }\right) _{\varepsilon }$ is a representative
of $u.$
\end{definition}

\begin{remark}
We have$\left( \widetilde{\partial }_{x_{j}}\right) _{h}^{0}u=u$ and for $%
\alpha \in \mathbb{Z}_{+}^{d},$
\begin{equation*}
\widetilde{\partial }_{h}^{\alpha }u=\left( \widetilde{\partial }%
_{x_{1}}\right) _{h}^{\alpha _{1}}\circ \left( \widetilde{\partial }%
_{x_{2}}\right) _{h}^{\alpha _{2}}\circ ...\circ \left( \widetilde{\partial }%
_{x_{d}}\right) _{h}^{\alpha _{d}}u\text{ }.
\end{equation*}
It is clear that $\widetilde{\partial }_{h}^{\alpha }u$ may be defined by
the representative $\left( \partial ^{\alpha }u_{\varepsilon }\ast \rho
_{h\left( \varepsilon \right) }^{\left[ \alpha \right] }\right)
_{\varepsilon }$, where $\rho _{h\left( \varepsilon \right) }^{\left[ \alpha %
\right] }=\rho _{h\left( \varepsilon \right) }\ast \rho _{h\left(
\varepsilon \right) }\ast ...\ast \rho _{h\left( \varepsilon \right) }$, the
convolution is taken $\left| \alpha \right| $ times. The notion of
regularized derivative is well defined and its class is independent of the
choice of the representative $\left( u\right) _{\varepsilon }$.
\end{remark}

In order to study the existence and uniqueness of Colombeau generalized
solutions of Cauchy problems with partial regularized derivatives, one
introduces the algebra of generalized functions suitable to this context.

We denote by $D_{L^{\infty }}\left( \overline{\Omega }\right) $ the algebra
of restrictions to $\overline{\Omega }$ of smooth functions defined on $%
\mathbb{R}^{d}$ with all derivatives bounded. With the same method of
construction of the simplified algebra of Colombeau, we define the
simplified algebra of global generalized functions, denoted $\mathcal{G}%
_{s,g}\left( \overline{\Omega }\right) ,$ by the quotient algebra
\begin{equation}
\mathcal{G}_{s,g}\left( \overline{\Omega }\right) =\frac{\mathcal{E}_{M,s,g}%
\left[ \overline{\Omega }\right] }{\mathcal{N}_{s,g}\left[ \overline{\Omega }%
\right] }\text{ \ ,}  \label{3.1}
\end{equation}
where
\begin{equation*}
\begin{array}{l}
\mathcal{E}_{M,s,g}\left[ \overline{\Omega }\right] =\left\{ \left(
u_{\varepsilon }\right) _{\varepsilon }\in \mathcal{E}_{s,g}\left[ \overline{%
\Omega }\right] :\forall \alpha \in \mathbb{Z}_{+}^{d},\exists p>0,\left\|
\partial ^{\alpha }u_{\varepsilon }\right\| _{L^{\infty }\left( \overline{%
\Omega }\right) }\leq O\left( \varepsilon ^{-p}\right) \right\} \\
\ \ \mathcal{N}_{s,g}\left[ \overline{\Omega }\right] =\left\{ \left(
u_{\varepsilon }\right) _{\varepsilon }\in \mathcal{E}_{s,g}\left[ \overline{%
\Omega }\right] :\forall \alpha \in \mathbb{Z}_{+}^{d},\forall q>0,\left\|
\partial ^{\alpha }u_{\varepsilon }\right\| _{L^{\infty }\left( \overline{%
\Omega }\right) }\leq O\left( \varepsilon ^{q}\right) \right\}
\end{array}
\end{equation*}
and $\mathcal{E}_{s,g}\left[ \overline{\Omega }\right] =\left( D_{L^{\infty
}}\left( \overline{\Omega }\right) \right) ^{I}$.

It is easy to see that $\mathcal{E}_{M,s,g}\left[ \overline{\Omega }\right] $
is a differential subalgebra of $\mathcal{E}_{s,g}\left[ \overline{\Omega }%
\right] $ and $\mathcal{N}_{s,g}\left[ \overline{\Omega }\right] $ is an
ideal of $\mathcal{E}_{M,s,g}\left[ \overline{\Omega }\right] $.

\begin{proposition}
\label{pro}Let $u\in D_{L^{\infty }}\left( \mathbb{R}^{d}\right) ,\alpha \in
\mathbb{Z}_{+}^{d}$ and $h\in \mathcal{H}$, if $h\left( \varepsilon \right)
=O\left( \varepsilon \right) ,\varepsilon \rightarrow 0$, then
\begin{equation}
\widetilde{\partial }_{h}^{\alpha }u=\partial ^{\alpha }u\text{ in }\mathcal{%
G}_{s,g}\left( \mathbb{R}^{d}\right) .  \label{3.2}
\end{equation}
\end{proposition}

\begin{remark}
In general if $u\in \mathcal{G}_{s}\left( \mathbb{R}^{d}\right) ,$ $%
\widetilde{\partial }_{h}^{\alpha }u\neq \partial ^{\alpha }u$ in $\mathcal{G%
}_{s,g}\left( \mathbb{R}^{d}\right) .$ For example, denote $H$ the Heaviside
function on $\mathbb{R}$, then
\begin{equation*}
\widetilde{\partial }_{h}H\neq H^{\prime }\text{ in }\mathcal{G}_{s,g}\left(
\mathbb{R}\right) .
\end{equation*}
\end{remark}

Let $T>0,$ $h\in \mathcal{H}$ and $W=\left[ -T,T\right] \times \mathbb{R}%
^{d} $ , the regularized derivative of an element $u$ of $\mathcal{G}%
_{s,g}\left( W\right) $ with respect to $x_{j}$ is defined as
\begin{equation*}
\left( \widetilde{\partial }_{x_{j}}\right) _{h}u=cl\left( \partial
_{x_{j}}u_{\varepsilon }\left( t,.\right) \ast \rho _{h\left( \varepsilon
\right) }\right) _{\varepsilon \in I},
\end{equation*}
where $\left( u_{\varepsilon }\right) _{\varepsilon \in I}$ is a
representative of $u$ and $\rho \in \mathcal{S}\left( \mathbb{R}^{d}\right) $
satisfies
\begin{equation}
\begin{array}{l}
\text{i) }\int \rho \left( x\right) dx=1 \\
\text{ii) }\int x^{\alpha }\rho \left( x\right) dx=0, \forall \alpha \in
\mathbb{N}^{d}
\end{array}
\label{dr1}
\end{equation}
Now we consider in $\mathcal{G}_{s,g}\left( W\right) $ the following linear
Cauchy problem
\begin{equation}
\left\{
\begin{array}{l}
\partial _{t}u+\sum\limits_{\left| \alpha \right| \leq m}a_{\alpha }%
\widetilde{\partial }_{h}^{\alpha }u=f \\
u\left( 0,x\right) =u_{0}(x)
\end{array}
\right.  \label{dr6}
\end{equation}
where $a_{\alpha }\in D_{L^{\infty }}\left( W\right) ,f\in \mathcal{G}%
_{s,g}\left( W\right) $ and $u_{0}\in \mathcal{G}_{s,g}\left( \mathbb{R}%
^{d}\right) .$

One of the main results of the paper \cite{col-he} is the following.

\begin{theorem}
\label{th-exis-unic}The linear Cauchy problem (\ref{dr6}) admits a global
unique solution $u\in \mathcal{G}_{s,g}\left( W\right) $ if there exists $%
p\in \mathbb{Z}_{+}$ such that
\begin{equation}
e^{C.h\left( \varepsilon \right) ^{-m}}=O\left( \varepsilon ^{-p}\right) ,
\label{dr7}
\end{equation}
where $c_{\alpha }=\left\| \partial ^{\alpha }\rho ^{\left[ \alpha \right]
}\right\| _{L^{1}\left( \mathbb{R}^{d}\right) }$ and $C=\sum\limits_{\left|
\alpha \right| \leq m}c_{\alpha }\left\| a_{\alpha }\right\| _{L^{\infty
}\left( W\right) }.$
\end{theorem}

\section{Non solvable differential operators}

The following differential operators
\begin{equation}
M=\frac{\partial }{\partial t}+ib\left( t\right) \frac{\partial }{\partial x}%
,  \label{M1}
\end{equation}
where $b\in C^{\infty }\left( \mathbb{R}\right) $ satisfies the condition
\begin{equation}
tb\left( t\right) >0,\forall t\in \mathbb{R}^{\ast },  \label{c.b}
\end{equation}
are called differential operators of Mizohata type. We know, see \cite{Ner},
that such operators $M$ are not locally solvable at the origin in the
framework of Schwartz distributions. A construction of a function $f\in
C_{0}^{\infty }\left( \mathbb{R}^{2}\right) $ such that there is no locally
distributional solution at the origin of the equation $Mu=f$ is given in
\cite{B-O}.

\begin{remark}
\bigskip It is well known that the operator (\ref{M1}) with the condition (%
\ref{c.b}) is reduced to the Mitzohata operator $\dfrac{\partial }{\partial t%
}+it\dfrac{\partial }{\partial x}$ if and only if $b\left( 0\right) =0$ and $%
b^{\prime }\left( 0\right) \neq 0$, see Tr\`{e}ves \cite{trev}. In our case
the function $b\left( t\right) $ may have a zero at the origin of infinite
order.
\end{remark}

In this section we give a necessary and sufficient condition for local
solvability of the equation
\begin{equation}
Mu\left( t,x\right) =f\left( t,x\right) ,  \label{C1}
\end{equation}
where $f\in C_{0}^{\infty }\left( \mathbb{R}^{2}\right) .$

Let $B\left( t\right) =\int_{0}^{t}b\left( s\right) ds$ and define the
function $Kf$ by
\begin{equation*}
Kf\left( x\right) =\int_{0}^{+\infty }\int_{-\infty }^{+\infty }e^{i\left(
x+iB\left( s\right) \right) \xi }\widehat{f}\left( s,\xi \right) dsd\xi ,
\end{equation*}
where $\widehat{f}\left( t,\xi \right) $ is the Fourier transform of $%
f\left( t,x\right) $ with respect to the variable $x.$

\begin{theorem}
\label{th-exis}The equation (\ref{C1}) admits a local distributional
solution at the origin of $\mathbb{R}^{2}$ if and only if the function $Kf$
is real analytic at the origin of $\mathbb{R}$.
\end{theorem}

\textbf{Proof} To solve the equation (\ref{C1}) we formally apply the
Fourier transformation with respect to $x,$ then
\begin{equation}
\frac{\partial \widehat{u}}{\partial t}\left( t,\xi \right) -b\left(
t\right) \xi \widehat{u}\left( t,\xi \right) =\widehat{f}\left( t,\xi
\right) ,  \label{4.4}
\end{equation}
hence
\begin{equation}
\widehat{u}\left( t,\xi \right) =\int_{t_{0}}^{t}e^{\left( B(t)-B(s)\right)
\xi }\widehat{f}\left( s,\xi \right) ds .  \label{4.5}
\end{equation}
To recover $u$ we must apply the inverse Fourier transformation to $\widehat{%
u}$, so the choice of $t_{0}$ is important. In (\ref{4.5}), we choose $t_{0}$
such that $\left( B(t)-B(s)\right) \xi \leq 0,$ $\forall \xi \in \mathbb{R}.$
By the condition (\ref{c.b}) the function $B$ is increasing for $t>0$ and
decreasing for $t<0$. \newline
For $\xi <0$, we choose $t_{0}=0$, and we define $u$ by
\begin{equation*}
\widehat{u}\left( t,\xi \right) =\int_{0}^{t}e^{\left( B\left( t\right)
-B\left( s\right) \right) \xi }\widehat{f}\left( s,\xi \right) ds .
\end{equation*}
For $\xi >0,$ we take
\begin{equation}
\widehat{u}\left( t,\xi \right) =\left\{
\begin{array}{c}
-\int_{t}^{+\infty }e^{\left( B\left( t\right) -B\left( s\right) \right) \xi
}\widehat{f}\left( s,\xi \right) ds,\text{ }t>0 \\
\int_{-\infty }^{t}e^{\left( B\left( t\right) -B\left( s\right) \right) \xi }%
\widehat{f}\left( s,\xi \right) ds,\text{ }t<0
\end{array}
\right.  \label{4.6}
\end{equation}
In this case the function $\widehat{u}$ admits a jump at $t=0$ given by
\begin{equation*}
\widehat{u}\left( +0,\xi \right) -\widehat{u}\left( -0,\xi \right)
=-\int_{-\infty }^{+\infty }e^{-B\left( s\right) \xi }\widehat{f}\left(
s,\xi \right) ds .
\end{equation*}
Consequently, we obtain in the distributional sense
\begin{equation}
\frac{\partial \widehat{u}}{\partial t}\left( t,\xi \right) -b\left(
t\right) \xi \widehat{u}\left( t,\xi \right) =\widehat{f}\left( t,\xi
\right) +\left[ \widehat{u}\left( 0_{+},\xi \right) -\widehat{u}\left(
0_{-},\xi \right) \right] \delta \left( t\right) \quad ,  \label{4.7}
\end{equation}
where $\delta $ is the Dirac measure at $0$. \newline
The inverse Fourier transform of (\ref{4.7}) with respect to $\xi $ gives
\begin{equation}
Mu\left( t,x\right) =f\left( t,x\right) -\delta \left( t\right) Kf\left(
x\right) .  \label{4.8}
\end{equation}
Let $H\left( t\right) $ be the Heaviside function$,$ then
\begin{equation*}
M\left( u\left( t,x\right) +H\left( t\right) Kf\left( x\right) \right)
=f\left( t,x\right) +ib\left( t\right) H\left( t\right) \left( Kf\left(
x\right) \right) ^{\prime }.
\end{equation*}
The term $ib\left( t\right) H\left( t\right) \left( Kf\left( x\right)
\right) ^{\prime }$ in the last equation is eliminated thanks to the
following function
\begin{equation*}
v\left( t,x\right) =\left\{
\begin{array}{c}
i\int_{0}^{t}b\left( s\right) \left( Kf\right) ^{\prime }\left( x-i\left(
B\left( t\right) -B\left( s\right) \right) \right) ds,\quad \text{ }t\geq 0
\\
0,\text{ \qquad \qquad \qquad \qquad \qquad \qquad \qquad \quad \quad \quad }%
t<0
\end{array}
\right.
\end{equation*}
which is well defined as $Kf\left( x\right) $ is assumed to be real analytic
at the origin. Therefore the function $v$ admits an holomorphic extension to
a neighborhood $\omega $ of the origin of $\mathbb{C}.$ Further the function
$v $ satisfies the equation $Mv(t,x)=ib(t)H(t)\left( Kf\right) ^{\prime
}(x). $ Define
\begin{equation}
w\left( t,x\right) =u\left( t,x\right) +H\left( t\right) Kf\left( x\right)
-v\left( t,x\right) ,  \label{4.9}
\end{equation}
then
\begin{equation*}
Mw\left( t,x\right) =f\left( t,x\right) \quad ,
\end{equation*}
i.e. $w\left( t,x\right) $ is a solution of the equation (\ref{C1}).

The constructed solution $w$ is of class $C^{\infty }$ in a neighborhood of
the origin. Indeed, we remark that if $t\neq 0$ the operator $M$ is elliptic
so $w$ is $C^{\infty }$ when $t\neq 0$. To show that $w$ is $C^{\infty }$ we
study the case $t=0$. We have

\begin{equation*}
w\left( t,x\right) =H\left( t\right) A\left( t,x\right) +H\left( -t\right)
B\left( t,x\right) ,
\end{equation*}
where
\begin{eqnarray*}
A\left( t,x\right) &=&\int_{-\infty }^{0}\int_{0}^{t}e^{\left( ix+B\left(
t\right) -B\left( s\right) \right) \xi }\widehat{f}\left( s,\xi \right)
dsd\xi -\int_{0}^{+\infty }\int_{t}^{+\infty }e^{\left( ix+B\left( t\right)
-B\left( s\right) \right) \xi }\widehat{f}\left( s,\xi \right) dsd\xi \\
&&+\int_{0}^{+\infty }\int_{-\infty }^{+\infty }e^{i\left( x+iB\left(
s\right) \right) \xi }\widehat{f}\left( s,\xi \right) dsd\xi
-i\int_{0}^{t}b\left( s\right) \left( Kf\right) ^{\prime }\left( x-i\left(
B\left( t\right) -B\left( s\right) \right) \right) ds
\end{eqnarray*}
and
\begin{equation*}
B\left( t,x\right) =\int_{-\infty }^{0}\int_{0}^{t}e^{\left( ix+B\left(
t\right) -B\left( s\right) \right) \xi }\widehat{f}\left( s,\xi \right)
dsd\xi +\int_{0}^{+\infty }\int_{-\infty }^{t}e^{\left( ix+B\left( t\right)
-B\left( s\right) \right) \xi }\widehat{f}\left( s,\xi \right) dsd\xi
\end{equation*}
It is clear that $w\left( t,x\right) $ is $C^{\infty }$ with respect to the
variable $x.$ Moreover we have
\begin{eqnarray*}
\partial _{t}^{j}\partial _{x}^{k}w\left( t,x\right) &=&H\left( t\right)
\partial _{t}^{j}\partial _{x}^{k}A\left( t,x\right) +H\left( -t\right)
\partial _{x}^{j}\partial _{t}^{k}B\left( t,x\right) \\
&&+\sum_{i=0}^{j-1}\delta ^{\left( i\right) }\left( t\right) \left( \partial
_{t}^{j-1-i}\partial _{x}^{k}w\left( 0_{+},x\right) -\partial
_{t}^{j-1-i}\partial _{x}^{k}w\left( 0_{-},x\right) \right) \\
&=&H\left( t\right) \partial _{t}^{j}\partial _{x}^{k}A\left( t,x\right)
+H\left( -t\right) \partial _{x}^{j}\partial _{t}^{k}B\left( t,x\right) \\
&&+\sum_{i=0}^{j-1}\delta ^{\left( i\right) }\left( t\right) \left( \partial
_{t}^{j-1-i}\partial _{x}^{k}A\left( 0_{+},x\right) -\partial
_{t}^{j-1-i}\partial _{x}^{k}B\left( 0_{-},x\right) \right) \text{ \ .}
\end{eqnarray*}
We also have
\begin{eqnarray*}
B\left( t,x\right) -A\left( t,x\right) &=&\int_{0}^{+\infty }\int_{-\infty
}^{+\infty }e^{\left( ix+B\left( t\right) -B\left( s\right) \right) \xi }%
\widehat{f}\left( s,\xi \right) dsd\xi -Kf\left( x\right) + \\
&&+i\int_{0}^{t}b\left( s\right) \left( Kf\right) ^{\prime }\left( x-i\left(
B\left( t\right) -B\left( s\right) \right) \right)ds,
\end{eqnarray*}
then $\forall l,k\in \mathbb{Z}_{+},\partial _{t}^{l}\partial _{x}^{k}\left(
B\left( t,x\right) -A\left( t,x\right) \right) $ is a finite sum of the
following terms
\begin{eqnarray*}
&&b^{\left( l_{1}\right) }\left( t\right) b^{l_{2}}\left( t\right) \left[
i^{k}\int_{0}^{+\infty }\int_{-\infty }^{+\infty }\xi ^{l_{3}+k}e^{\left(
ix+B\left( t\right) -B\left( s\right) \right) \xi }\widehat{f}\left( s,\xi
\right) dsd\xi \right.  \notag \\
&&\left. -\left( -i\right) ^{l_{3}}\left( Kf^{\left( l_{3}+k\right) }\left(
x\right) -i\int_{0}^{t}b\left( s\right) \left( Kf\right) ^{\left(
l_{3}+1\right) }\left( x-i\left( B\left( t\right) -B\left( s\right) \right)
\right) ds\right) \right]
\end{eqnarray*}
where $l_{1},l_{2}$ and $l_{3}$ depend only on $l$. It is clair that these
terms equal all zero when $t=0$, then
\begin{equation*}
\partial _{t}^{j}\partial _{x}^{k}w\left( t,x\right) =H\left( t\right)
\partial _{t}^{j}\partial _{x}^{k}A\left( t,x\right) +H\left( -t\right)
\partial _{x}^{j}\partial _{t}^{k}B\left( t,x\right)
\end{equation*}
which give $w\in C^{\infty }$.

The proof of the necessity of the analyticity of $Kf$. Let us suppose that
there is $u\in C^{1}$ such that $Mu=f$ in a neighborhood $\Omega $ of the
origin and let $\chi \in C_{0}^{\infty }\left( \mathbb{R}^{2}\right) ,\chi
\equiv 1$ in a neighborhood of the origin and $supp\chi \subset \Omega .$ So
$\left( \chi f\right) $ satisfies locally $Mu=\chi f$\ and therefore the
function $Kf(x)$ can be written in the form
\begin{equation*}
Kf(x)=\underset{\varepsilon \rightarrow 0}{\lim }\int_{-\infty }^{+\infty
}\int_{-\infty }^{+\infty }\int_{0}^{+\infty }e^{i\left( x-y+iB\left(
s\right) \right) \xi -\varepsilon \xi ^{2}}\left( \chi f\right) (s,y)\frac{%
d\xi }{2\pi }dyds\text{ .}
\end{equation*}
Consider the integral
\begin{equation*}
K_{\varepsilon }f(x)=\int_{-\infty }^{+\infty }\int_{-\infty }^{+\infty
}\int_{0}^{+\infty }e^{i\left( x-y+iB\left( s\right) \right) \xi
-\varepsilon \xi ^{2}}\left( \chi f\right) (s,y)\frac{d\xi }{2\pi }dyds\ ,
\end{equation*}
then
\begin{eqnarray}
K_{\varepsilon }f(x) &=&\int_{-\infty }^{+\infty }\int_{-\infty }^{+\infty
}\int_{0}^{+\infty }e^{i\left( x-y+iB\left( s\right) \right) \xi
-\varepsilon \xi ^{2}}M(\chi u)\left( s,y\right) \frac{d\xi }{2\pi }dyds
\notag \\
&&-\int_{-\infty }^{+\infty }\int_{-\infty }^{+\infty }\int_{0}^{+\infty
}e^{i\left( x-y+iB\left( s\right) \right) \xi -\varepsilon \xi ^{2}}u(M\chi
)\left( s,y\right) \frac{d\xi }{2\pi }dyds\text{ .}  \label{4.11}
\end{eqnarray}
An integration by parts of the first term of the second member gives
\begin{equation}
\int_{-\infty }^{+\infty }\int_{-\infty }^{+\infty }\int_{0}^{+\infty
}e^{i\left( x-y+iB\left( s\right) \right) \xi -\varepsilon \xi ^{2}}M(\chi
u)\left( s,y\right) \frac{d\xi }{2\pi }dyds=0,
\end{equation}
hence
\begin{equation*}
K_{\varepsilon }f(x)=-\int_{-\infty }^{+\infty }\int_{-\infty }^{+\infty
}\int_{0}^{+\infty }e^{i\left( x-y+iB\left( s\right) \right) \xi
-\varepsilon \xi ^{2}}\frac{d\xi }{2\pi }u(s,y)M\chi \left( s,y\right) dyds%
\text{ .}
\end{equation*}
Consider now the deformation of the path of integration with respect to $\xi
$ in $K_{\varepsilon }f(x)$ by taking the contour $\Gamma $ defined by
\begin{equation*}
\zeta =\rho \left( 1+\frac{i}{2}\frac{x-y}{\left| x-y\right| }\right) \quad ,%
\text{ \ }\rho >0.
\end{equation*}
Hence, for any $\varepsilon >0$ fixed, we have
\begin{equation}
K_{\varepsilon }f(x)=-\int_{-\infty }^{+\infty }\int_{-\infty }^{+\infty }%
\underset{\Gamma }{\int }e^{i\left( x-y+iB\left( s\right) \right) \zeta
-\varepsilon \zeta ^{2}}\frac{d\zeta }{2\pi }u(s,y)M\chi \left( s,y\right)
dyds .  \label{4.10}
\end{equation}
The function $K_{\varepsilon }f\left( x\right) $ is analytic in $x$ for each
fixed $\varepsilon .$ It remains to show that $\lim\limits_{\varepsilon
\rightarrow 0}K_{\varepsilon }f\left( x\right) $ is analytic at the origin.
For this need, we have to estimate uniformly the expression $K_{\varepsilon
}f\left( x\right) .$ Since $M\chi =0$ in a neighborhood of the origin, as $%
\chi \equiv 1$ in this neighborhood, so in (\ref{4.10}) the integral with
respect to $s$ and $y$ is taken outside a rectangle, i.e. either $\left|
s\right| >c_{1}$ or $\left| y\right| >c_{2}.$ Consequently, $B\left(
s\right) >c_{3}$ or $\left| x-y\right| >c_{4},$ where $\left( c_{j}\right)
_{j=1}^{4}$ are positive constants not depending on $s,$ $y$ $,$ $x$. Then
\begin{equation*}
Im\left( \left( x-y+iB(s)\right) \zeta -i\varepsilon \zeta ^{2}\right) \geq
c\rho +\frac{3}{4}\varepsilon \rho ^{2},
\end{equation*}
from this estimate we conclude that $\underset{\varepsilon \rightarrow 0}{%
\lim }K_{\varepsilon }f(x)$ is analytic with respect to $x$ in a
neighborhood of the origin.

Now suppose that $u\in D^{\prime }\backslash C^{1}$ and $u$ is a solution of
$Mu=f$, then $u$ is a $C^{\infty }$ function of $t$ with values in $%
D^{\prime }\left( \mathbb{R}\right) ,$ see theorem 4.4.8 \cite{Hor}.\ By the
local structure of distribution, we can assume that, there exists a function
$v\in C^{1}\left( \mathbb{R}^{2}\right) $ such that $u=\partial _{x}^{N}v$.
As $M\left( \partial _{x}^{N}v\right) =\partial _{x}^{N}\left( Mv\right) ,$
we may substitute $\partial _{x}^{N}v$ for $u$ in (\ref{4.11}) and proceed
in a same way to obtain the general result.\hfill $\square $

\section{Differential operators of Mizohata type in $\mathcal{G}_{s,g}$}

Consider in $\mathcal{G}_{s,g}\left( W\right) ,W=\left[ -T,T\right] \times
\mathbb{R},$ the following equation
\begin{equation}
\partial _{t}U+ib\left( t\right) \widetilde{\partial _{x}}_{h}U=f\text{ in }%
\mathcal{G}_{s,g}\left( W\right) ,  \label{E1}
\end{equation}
where $f\in C_{0}^{\infty }\left( \mathbb{R}^{2}\right) $, $h\in \mathcal{H}$%
, $b\in D_{L^{\infty }}\left( \mathbb{R}\right) $ and
\begin{equation*}
tb\left( t\right) >0,t\in \mathbb{R}^{\ast }\text{ ,}
\end{equation*}
Then we have the following result.

\begin{theorem}
\label{thf1}Let $U=cl\left( u_{\varepsilon }\right) _{\varepsilon }\in
\mathcal{G}_{s,g}\left( W\right) $ be a solution of (\ref{E1}) which is
locally associated to a distribution $v\in D^{\prime }\left( \omega \right) $
at the origin, then the function $Kf$ is analytic in a neighborhood of the
origin of $\mathbb{R}$.
\end{theorem}

\textbf{Proof} Let us suppose that a solution $U=cl\left( u_{\varepsilon
}\right) _{\varepsilon }$ of (\ref{E1}) is locally associated to a
distribution $v$ at the origin, then there is a neighborhood $\widetilde{%
\omega }$ of the origin such that, $\forall \phi \in C_{0}^{\infty }\left(
\widetilde{\omega }\right) $ ,
\begin{equation*}
\lim_{\varepsilon \rightarrow 0}\int u_{\varepsilon }\left( t,x\right) \phi
\left( t,x\right) =\left\langle v,\phi \right\rangle .
\end{equation*}
We have

\begin{equation*}
\left( \partial _{t}u_{\varepsilon }\right) _{\varepsilon }+ib\left(
t\right) \left( \partial _{x}u_{\varepsilon }\ast \rho _{h\left( \varepsilon
\right) }\right) _{\varepsilon }-f\in \mathcal{N}_{s,g}\left[ W\right] ,
\end{equation*}
where the convolution takes place in the $x$-variable at fixed $t$, so
\begin{equation*}
\underset{\varepsilon \rightarrow 0}{\lim }\left( \left( \partial
_{t}u_{\varepsilon }\right) _{\varepsilon }+ib\left( t\right) \left(
\partial _{x}u_{\varepsilon }\ast \rho _{h\left( \varepsilon \right)
}\right) _{\varepsilon }-f\right) =0\text{ in }D^{\prime }\left( \widetilde{%
\omega }\right) .
\end{equation*}
Since $\forall h\in \mathcal{H},$ the sequence $(\rho _{h\left( \varepsilon
\right)})_ \varepsilon $ converges to the Dirac measure, as $\varepsilon
\rightarrow 0$, then $\forall \phi \in C_{0}^{\infty }\left( \widetilde{%
\omega }\right) $ we have
\begin{equation*}
\lim_{\varepsilon \rightarrow 0}\int \left( \partial _{x}u_{\varepsilon
}\ast \rho _{h\left( \varepsilon \right) }\right) _{\varepsilon }\left(
t,x\right) \phi \left( t,x\right) =\left\langle \partial _{x}v,\phi
\right\rangle ,
\end{equation*}
hence $\partial _{t}v+ib\left( t\right) \partial _{x}v=f$ in $D^{\prime
}\left( \widetilde{\omega }\right) ,$ i.e. $v$ is a solution of the equation
\begin{equation}
Mu=f\text{ in }D^{\prime }\left( \widetilde{\omega }\right) ,  \label{E2}
\end{equation}
consequently, by theorem \ref{th-exis}, the function $Kf$ is analytic in
neighborhood of the origin of $\mathbb{R}$.\hfill $\square $

\begin{theorem}
\label{thf2}If $Kf$ is analytic in a neighborhood of the origin of $\mathbb{R%
}$, and the distributional solution $v\in D^{\prime }\left( \omega
\right)$ of (\ref{C1}) satisfies $(\widetilde{\partial
}_{x})_{h}v=\partial_{x} v$, then (\ref{E1}) admits a solution
$U=cl\left( u_{\varepsilon }\right) _{\varepsilon }\in
\mathcal{G}_{s,g}\left( W\right) $ which is locally associated to
$v$ at the origin.
\end{theorem}

\textbf{Proof} Let us suppose that the function $Kf$ is analytic in a
neighborhood of the origin$,$ and let $v\in D^{\prime }\left( \omega \right)
$ such that
\begin{equation*}
Mv=f\text{ in }D^{\prime }\left( \omega \right) .
\end{equation*}
Moreover, see the proof of the theorem \ref{th-exis},$\ v$ is of class $%
C^{\infty }$ in $\omega .$ Let $\Omega \Subset \omega ,$ then $v\in
D_{L^{\infty }}\left( \overline{\Omega }\right) $ and by the fact that $%
\left( \widetilde{\partial _{x}}\right) _{h}v=\partial _{x}v$ we have
\begin{equation}
\partial _{t}v+ib\left( t\right) \left( \widetilde{\partial _{x}}\right)
_{h}v=f\text{ in }\mathcal{G}_{s,g}\left( \overline{\Omega }\right) .
\label{5.3}
\end{equation}
Let $U$ be a solution of the equation (\ref{E1}) in $\mathcal{G}_{s,g}\left(
W\right) $ with the initial data given by $U\left( 0,x\right) =\left[
v\left( 0,x\right) \right] ,$ then in $\mathcal{G}_{s,g}\left( \overline{%
\Omega }\right) $ we have
\begin{equation*}
\left\{
\begin{array}{l}
\partial _{t}\left( U-v\right) +ib\left( t\right) \left( \widetilde{\partial
_{x}}\right) _{h}\left( U-v\right) =0, \\
\left( U-v\right) \left( 0,x\right) =0.
\end{array}
\right. \text{ }
\end{equation*}
The uniqueness of the generalized solution in the theorem \ref{th-exis-unic}
gives $U-v=0$ in $\mathcal{G}_{s,g}\left( \overline{\Omega }\right) $, hence
$U\approx _{0}v.$\hfill $\square $

\begin{remark}
The condition $\left( \widetilde{\partial _{x}}\right)
_{h}v=\partial _{x}v$ may be replaced by the condition $h\left(
\varepsilon \right) =O\left( \varepsilon \right) ,\varepsilon
\rightarrow 0$.
\end{remark}

For the Mizohata equations under consideration the new generalized
solutions from the method in \cite{col-he} can be associated with
distributions only in the case these equations are solvable in the
sense of distributions theory. This follows at once from theorems
\ref{th-exis}, \ref{thf1} and \ref {thf2}.

\end{document}